\documentclass[12pt]{article}
\usepackage[english]{babel}
\usepackage[utf8]{inputenc}
\usepackage{amsmath}
\usepackage{graphicx}
\usepackage{color}
\usepackage{setspace}
\usepackage{enumerate}
\usepackage{amssymb}
\usepackage{blindtext}
\numberwithin{equation}{section}
\newcommand{\Gve}{\varepsilon}
\begin{document}
\title{On ideal dynamic climbing ropes}
\author{D. Harutyunyan, G.W. Milton,  T.J. Dick and J. Boyer\\
\textit{Department of Mathematics, The University of Utah}}

\date{\today}
\maketitle

\begin{abstract}
We consider the rope climber fall problem in two different settings. The simplest formulation of the problem is when the climber falls from a given altitude and is attached to one end of the rope while the other end of the rope is attached to the rock at a given height. The problem is then finding the properties of the rope for which the peak force felt by the climber during the fall is minimal. The second problem of our consideration is again minimizing the same quantity in the presence of a carabiner. We will call such ropes \textit{mathematically ideal.}
Given the height of the carabiner, the initial height and the mass of the climber, the length of the unstretched rope, and the distance between the belayer and the carabineer, we find the optimal (in the sense of minimized the peak force to a given elongation) dynamic rope in the framework of nonlinear elasticity. Wires of shape memory materials have some of the desired features of the tension-strain relation of a mathematically ideal dynamic rope, namely a plateau in the tension over a range of strains. With a suitable hysteresis loop, they also
absorb essentially all the energy from the fall, thus making them an ideal rope in this sense too.
\end{abstract}


\section{Introduction}
\label{sec1}
Climber fall is a central problem in rock climbing; and an important factor, which this paper addresses, is minimizing
the peak force felt by the climber as he/she falls, allowing a maximum elongation of the rope.
Minimizing the peak force is important also for decreasing the likelihood that the anchor will be dislodged, and for minimizing
the stress on the rope and on the bolts, carabiners and all kind of other protections, thus increasing their lifetime.
Some comprehensive analyses of the maximal forces and rope elongation and the design of optimal ropes, optimality meant in various senses, already exist. In [\ref{bib:Leu1.},\ref{bib:Leu2.}] Leuth\"auser studies the above general problem in the setting of viscoelasticity. When a climber falls, some of the energy
is converted to heat during the fall, and it is assumed in both [\ref{bib:Leu1.},\ref{bib:Leu2.}], that the
coefficient of conversion, i.e, the proportion of the energy converted to heat and the total energy,
is equal to $0.5.$ In [\ref{bib:Sporri.}] Sp\"orri carries out a numerical analysis for the forces
 acting on the climber, and the resultant motion of the climber,
assuming Pavier's 3-parameter viscoelastic rope model [\ref{bib:Pav1.}].

There are other problems, too, that concern climber fall.
One is the durability of the rope and a safe lower  bound on the number of the falls the rope can
handle due to tear during any fall, which has been studied numerically by Bedogni and Manes [\ref{bib:Bed.Man.}], and in combination
with experiment by Pavier [\ref{bib:Pav1.}]. A second problem concerns the stiffness or stretchiness of the rope
(which controls the total extension of the rope during a fall),
and the maximal handling force of bolts and carabiners, see [\ref{bib:Att.},\ref{bib:Bla.Cus.Gra.Oka}].
A third problem is mostly medical and concerns the most probable injuries of the climber, [\ref{bib:Pai.Fio.Hou.}].

There are two kinds of ropes used in rope climbing: static and dynamic.
A dynamic rope stretches and is the rope connecting the climber to the belayer or anchor, while a
static rope has little stretch and is often used in ascending using jumars, hauling, and for
rappelling, equivalently called abseiling. In addition, static ropes or webbing
are used to affix via carabiners the dynamic rope to the climbing wall.
The problem we consider is to find and minimize the maximal force acting on the climber during a fall over a fixed distance. This peak force
is the maximal value of the tension of the dynamic rope at the point attached to the climber.
In other words, the first goal is to identify the properties of a "mathematically ideal dynamic rope" that is designed so the peak force on the climber is minimized during a fall with a specified elongation of the rope. Of course many factors influence whether a dynamic rope
is ideal in practice (such as durability, cost, weight, knotability and redundancy so
breaks of single filaments do not affect the overall strength), and for this reason we use the term ``mathematically ideal'' to mean
ideal in the restricted sense of minimizing the peak force for a given total elongation.
A second goal is to obtain a rope which absorbs essentially all the energy of the fall, so that the climber does not
rebound after the fall in contrast to the oscillations that bungy cords produce during bungy jumping.
To achieve the first goal, we begin by assuming the rope is purely elastic, rather than
viscoelastic. Later we consider such ropes with suitable hysteresis loops which absorb all the energy from the fall, and thus which are ideal
from the viewpoint of both goals.

The question is then how the nonlinear elastic response of the rope should be tailored to minimize the maximal force acting on the climber during a fall.
Intuitively, it makes sense that the rope should be designed to provide a constant braking force on the falling climber. While it may be the case that
this result is known, we were unable to find a reference which addressed this condition. So a mathematical proof of this fact is provided.
We draw attention to the fact that shape memory material wires have the desired characteristic of a plateau in the stress as the strain is varied.
Furthermore, they can exhibit large hysteresis loops, and we observe that this feature is exactly what is needed to achieve the second goal,
i.e., to absorb all the energy of the fall without oscillation.

A similar problem has been studied by Reali and Stefanini in [\ref{bib:Rea.Ste.}] in the setting of linear elasticity, but we believe nonlinear elasticity is more appropriate to the analysis of dynamic climbing ropes. We also remark that the problem of minimizing the maximal force during
deceleration over a fixed distance is also appropriate to pilots landing on an aircraft carrier, but there the desired response can be controlled
through the hydraulic dampers attached to the braking cable.


\section{Rope without a carabiner}
\label{sec2}
\subsection{Notation and problem setting}
\label{subsec2.1}
\begin{figure}[h]
\centering
\includegraphics[width=0.9\textwidth]{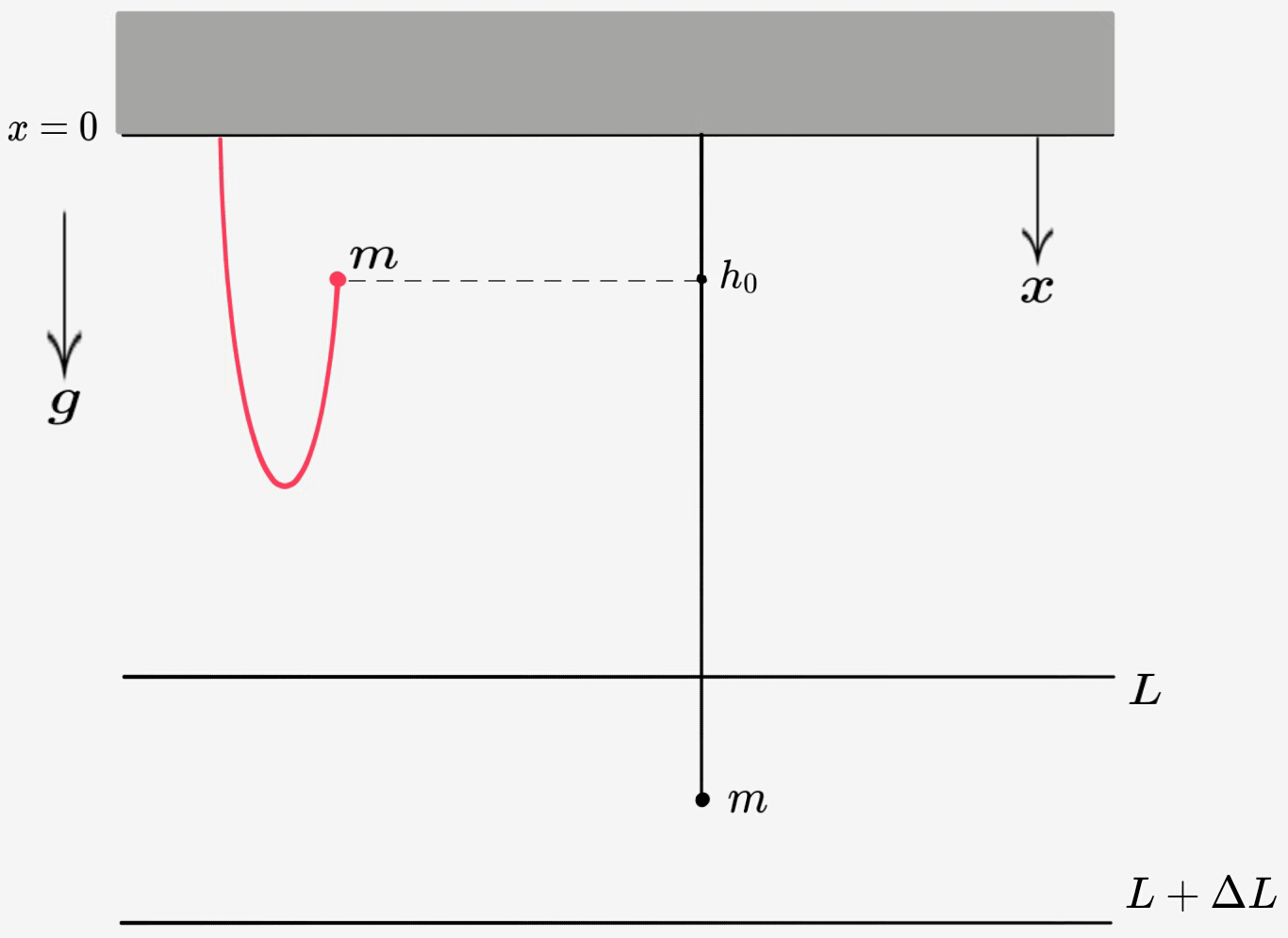}
\caption{\label{Fig1}Before and after the fall. The figure is schematic, in that the climber before the fall, as illustrated on the left, should be
directly below or above the point where the rope is attached.}
\end{figure}

In this section, we set up the notation for the climber fall problem in its simplest formulation. We direct the $x$ axis towards the direction of gravity: see Figure~\ref{Fig1}. We shall choose the gravitational potential to be zero at the position $x=0$ for convenience. Let us mention, that
the variable $x$ does not show the position of the climber, but rather it is the coordinate variable used in the undeformed configuration.
There will be a climber affixed to a rock wall by means of a rope and the climber will fall vertically downwards from a point directly above or below where the rope is attached to the wall, so that the problem is one-dimensional. The given rope of length $L$ will be attached to the rock at the point $x=0$ and the climber of mass $m$ will be attached
 to the other end at some arbitrary height $x=h_{0}$, where \ $|h_{0}|\leq|L|$. Note that $h_0$ denotes the distance below where the rope is anchored, and is negative if the climber is above the anchoring point. The rope will be assumed to be massless.
We assume that the climber is allowed to fall without hitting the ground up to the point $L+\Delta L,$ i.e., the maximal admissible stretching of the rope
is $\Delta L.$\\
\textit{We are interested in finding the properties of the rope so that given the initial data $L, \Delta L, m$ and $h_0$ it minimizes the maximal force felt by the climber during a fall. Such a rope will be called \textbf{mathematically ideal}.}\\
 We reformulate the problem as follows: instead of considering what happens before the rope becomes taut, we can set our clock so $t=0$  marks the instant in time when the rope becomes taut; and at that time the gravitational potential energy
 $mg(L-h_0)$ lost by the climber in falling a distance $L-h_0$ will have been converted into kinetic energy so that
the climber will have a velocity $v_0=\sqrt{2g(L-h_0)}$ in the direction of increasing $x$ (downwards) at time $t=0$. Thus the deformation of the rope begins at time $t=0$. We denote by $y(x,t)$ the position of a point on the rope at time $t$ that was
initially located at $x$ at time $t=0;$ thus a point on the rope at $x$ at $t=0$ gets displaced by a distance $u(x,t)=y(x,t)-x$
 (here $y(x,t)$ will be assumed to be continuous and differentiable in $x$).

Once the rope begins to deform, its deformation is assumed to be described by nonlinear elasticity theory (ignoring
viscosity).
\begin{figure}[h]
\centering
\includegraphics[width=0.9\textwidth]{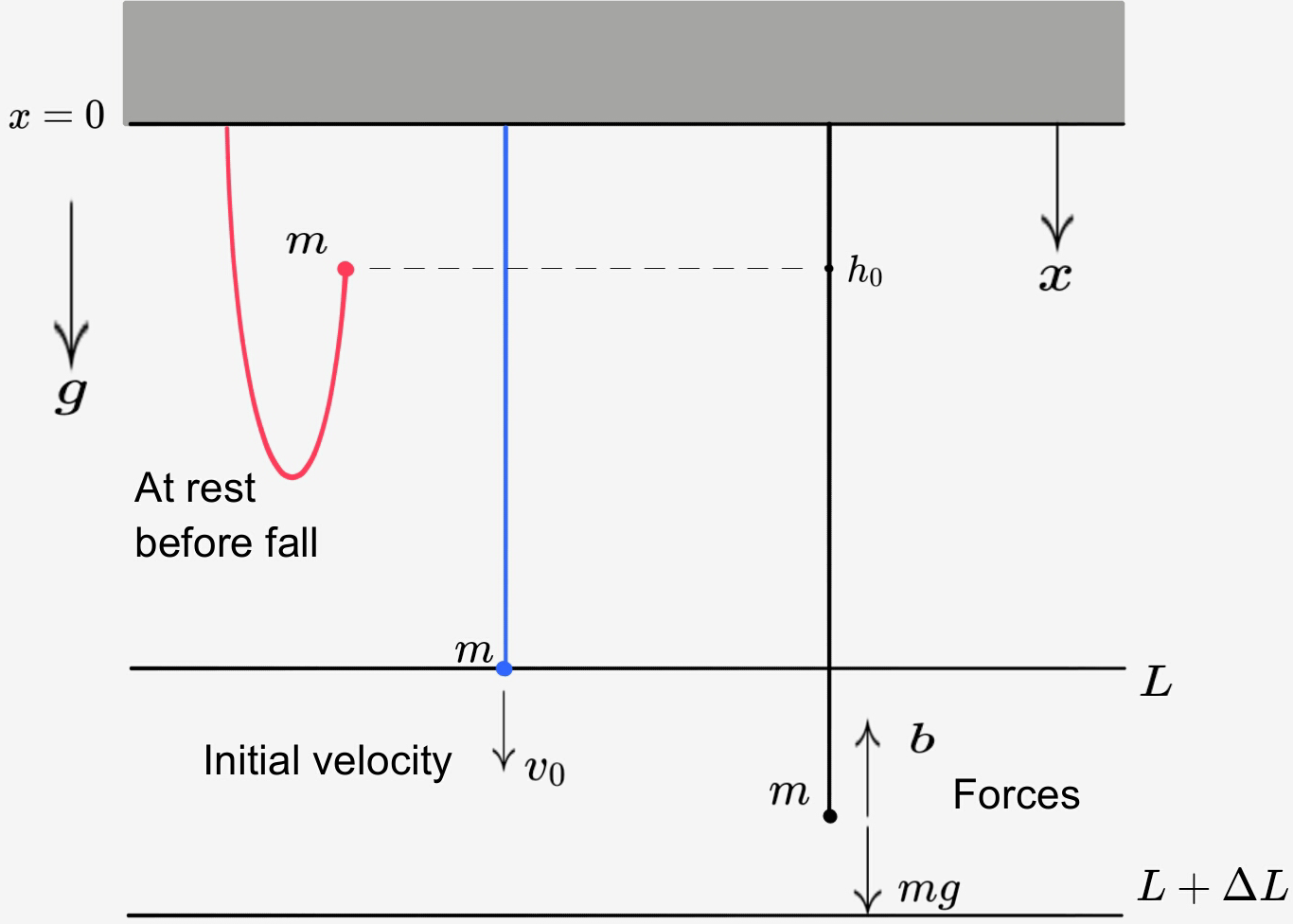}
\caption{\label{Fig2}The reformulated problem.}
\end{figure}
 The elastic properties of the rope are given by a function $W$, representing the elastic energy density, of the one-dimensional strain
$$ \Gve(x)=\frac{\partial u(x,t)}{\partial x}=\frac{\partial y(x,t)}{\partial x}-1. $$
We assume the elastic energy does not depend on higher order derivatives of the deformation $y(x,t)$, that there is no air resistance, and (to begin with, as it will be something we will reconsider later in section 4) that there is
no energy absorption or dissipation during a fall. Thus, a fall will be a periodic phenomenon as the climber oscillates between the heights $x=h_0$ and $x=L+\triangle L.$ We denote the time the climber reaches the critical point $x=L+\triangle L$ by $t=T.$ We denote furthermore by $E_{el}$ the elastic energy of the rope, and by $E_{total}$ the total energy of the system.
The total elastic energy for the rope, treated as a one dimensional body, is then:
\begin{equation}
\label{2.2}
E_{el}(t)=\int \limits_{0}^{L} W\left(\Gve(x)\right) \, dx. \quad
\end{equation}
In other words, the total elastic energy is the sum (integral) of the elastic energies associated with the stretching of each rope segment. Note that we do {\it not} assume that $W$ is a quadratic function of $\Gve$ and hence it is not appropriate to associate an elastic modulus to the rope: the elastic modulus is an appropriate descriptor when the
tension in the rope is proportional to the overall strain (and it is the elastic modulus which gives this constant of proportionality). However, in our {\it mathematically ideal} rope, we will see it is
best to have a rope where the tension is independent of the strain.

From the theory of nonlinear elasticity, e.g, [\ref{bib:Gur.}] the elastic energy density $W$ has to satisfy the following properties:
\begin{itemize}
\item[(\textbf{P1})] $W(\Gve)\geq 0$ \ for all \ $\Gve\in\mathbb R$ (deformations store elastic energy)
\item[(\textbf{P2})] $W(0)=0$ (no deformation-no energy)
\item[(\textbf{P3})] $W(\Gve)$ is a quasiconvex function of $\Gve+1$, which is equivalent to being a convex function of $\Gve$ in the one-dimensional case under consideration. Recall, that a function $f(x)\colon[a,b]\to\mathbb R$ is convex, if it satisfies the inequality
    $$f(\lambda x+(1-\lambda)y)\leq \lambda f(x)+(1-\lambda)f(y),$$
   for all $x,y\in[a,b],\lambda\in [0,1].$
Geometrically, convexity means that the entire graph of the function is above the tangent line at any point $x\in[a,b].$ Physically, convexity is important for stability: if  $W$ is not convex,
then microscopic oscillations in the deformation of the rope are energetically favorable, and macroscopically the behavior of the rope will be given
by an energy function which is again convex, being the convexification of the local nonconvex energy function.  (Such oscillatons do occur in shape memory wires, a point which we will return to later.)
\end{itemize}
 We aim to find under which conditions on $W$ the rope is \textit{mathematically ideal.} Our goal is to mitigate the force on the climber during the fall. In the more general version of the problem, the only new feature will be that the rope also passes through one carabiner, attached to the rock, that has a friction coefficient $k$.  The single carabiner case will be sufficient to make a generalization with an arbitrary integer number of carabiners, with frictional coefficients $k_{i}$.\\

\subsection{Optimal bounds via energy conservation}
\label{subsec2.2}
Denote by $b(t)$ the rope tension at the point $x=L$ at time $t,$ and by $a(x,t)$ the acceleration of the rope point $x$ at time $t,$ for $x\in[0,L]$ and $t\in\mathbb R.$ We have to solve the minimization problem
\begin{equation}
\label{2.3}
\min(\max_{t\geq 0}|b(t)|).
\end{equation}
It is clear, as the motion of the climber is periodic, we can consider the minimization problem on the restricted time interval $[0,T]$ rather than the whole
time interval $[0,\infty).$
It is evident that the tension $b(t)$ imposes an upwards force on the climber for all $t\in[0,T].$ The other force acting on the climber is the gravitation force $mg$ pointing downwards, thus being positive. So, the net force on the climber is $mg-b(t)$ and by Newton's second law we have that the mass times the acceleration downwards is given by
\begin{equation}
\label{2.4}
    ma(L,t)=mg-b(t).
\end{equation}
Next, we prove the key inequality
\begin{equation}
\label{2.5}
\max_{t\in[0,T]}|b(t)| \geq \frac{mg(L+\Delta L-h_0)}{\Delta L}.
\end{equation}
It is actually a direct consequence of energy conservation. Namely, as is well-known, we have that the work done by the climber in the time interval $[0,T]$ equals on one hand $\int_{L}^{L+\Delta L}ma(L,t)dx,$ and on the other hand it is the change in kinetic energy, i.e.,
$-\frac{mv_0^2}{2}.$ Therefore we get
\begin{equation}
\label{2.6}
\int_{L}^{L+\Delta L}ma(L,t)dx=-\frac{mv_0^2}{2}.
\end{equation}
Integrating (\ref{2.4}) in $x$ over the interval $[L,L+\Delta L],$ we obtain
$$ -\int_{L}^{L+\Delta L}ma(L,t)dx=\int_{L}^{L+\Delta L}b(t)dx-mg\Delta L,$$
and thus taking into account (\ref{2.6}) and the formula for the initial velocity, $v_0=\sqrt{2g(L-h_0)}$, we get
\begin{equation}
\label{2.7}
\int_{L}^{L+\Delta L}b(t)dx=mg\Delta L+\frac{mv_0^2}{2}=mg(L+\Delta L-h_0).
\end{equation}
Finally applying the inequality
$$\left|\int_{L}^{L+\Delta L}b(t)dx\right|\leq \Delta L\max_{t\in[0,T]}|b(t)|$$
to (\ref{2.7}), we arrive at (\ref{2.5}). Observe, that equality in (\ref{2.5}) holds if and only if
the tension $b(t)$ and thus the acceleration $a(L,t)$ is constant in the interval $[0,T],$ which is the scenario that
a \textit{mathematically ideal rope} must develop. Thus in this case we get
\begin{align}
\label{2.8}
a(L,t)=a_0&\equiv \frac{g(h_0-L)}{\Delta L},\\ \nonumber
b(t)=b_0&\equiv\frac{mg(L+\Delta L-h_0)}{\Delta L}.
\end{align}
This \textit{mathematically ideal rope} causes the force where the rope is attached to the climber to step from a zero force
to a constant value as the rope becomes taut. In practice, one would want a more gradual transition to allow the climber's body
time to respond to the force. Certainly the effects of the force cannot propagate through the climber's body faster than the speed of elastic waves; but
more significantly, if we consider the climber's body as a viscoelastic object, then one would not want the transition time to be faster than
the typical viscoelastic relaxation time. The undesireable instantaneous step function in the force will be mollified if we replace the \textit{mathematically ideal rope} by an approximation to it, such as the shape memory material rope, suggested later, which has a transition region before the stress plateau.

\subsection{The Optimal Elastic Energy Density Function}
\label{subsec2.3}
In this section, we find a formula for the elastic energy density function $W$ associated with a {\it mathematically ideal rope}.
Under our assumption that the rope is massless, equilibrium of forces implies the tension $b$ in the rope must be constant along the rope. Also if the rope is
optimal in the sense that the maximal possible force felt by the climber is minimized, then $b$ must be independent of the strain and given by (\ref{2.8}). Suppose a small segment of rope extending from $x=x_0$ to $x=x_0+\ell$ in the undeformed state at $t=0$, gets extended under the deformation at some time $t>0$ less than $T$ to the length $\ell+\delta\ell$. The work done on the rope $b_0\delta\ell$ (being the force times the distance)
must go into the elastic energy stored in this rope segment, implying
\begin{equation}
\label{2.14.5}
b_0\delta\ell=\int \limits_{x_0}^{x_0+\delta\ell} W\left(\frac{\partial y(x,t)}{\partial x}-1\right) \, dx.
\end{equation}
Using the fact that the energy density is zero in the undeformed state at $t=0$ when $y(x,t)=x$, this will hold if
\begin{equation}
\label{2.15}
W(\Gve)=b_0\Gve=\frac{mg(L+\Delta L-h_0)}{\Delta L}\Gve.
\end{equation}
This defines $W$ for values $\Gve\geq 0$. (Equivalently, one could use the fact that the local tension $b=b_0$, which is the one-dimensional stress,
is $\partial W(\Gve)/\partial \Gve$ to deduce that $W(\Gve)$ necessarily has the form (\ref{2.15})).
Moreover, it is clear that $W$ satisfies the properties (P1) and (P2) for $\Gve\geq 0.$ The case $\Gve<0$ corresponds
to rope compression, which is assumed to require no energy, thus we must take $W(\Gve)=0$ for $\Gve<0.$
In conclusion, we get a final formula for the energy density function:
\begin{eqnarray}
\label{2.16}
W(\Gve) & = & \frac{mg(L+\Delta L-h_0)}{\Delta L}\Gve ~~\text{if} \ \ \ \Gve\geq 0, \\ \nonumber
        & = & 0 ~ \text{if} \ \ \ \Gve<0,
\end{eqnarray}
which clearly satisfies all properties (P1-P3).

Note however that $W(\Gve)$ is not a strictly convex function of $\Gve$. If it was strictly convex, we could use
Jensen's inequality,
\begin{equation}
\label{2.16.1}
 \frac{1}{L}\int \limits_{0}^{L} W\left(\frac{\partial y(x,t)}{\partial x}-1\right) \, dx\geq W\left(\frac{y(L,t)}{L}-1\right),
\end{equation}
with equality holding only when $y(x,t)$ is independent of $x$, to deduce that the deformation of the rope which minimizes the elastic energy is necessarily
homogeneous\footnote{A function $f(x,t)\colon\mathbb R^2\to\mathbb R$ is homogeneous in $x$ for each $t\in\mathbb R,$ if $f(\lambda x,t)=\lambda f(x,t)$ for all $\lambda,x,t\in\mathbb R.$ In other words, its derivative in the $x$ variable does not depend on $x.$}
at each time $t$. Recall, that Jensen's inequality asserts the following: \textit{If a function $f\colon\mathbb R\to\mathbb R$ is convex, then
for any interval $[a,b]\subset\mathbb R$ and any continuous function $g\colon[a,b]\to\mathbb R$ the inequality holds:}
\begin{equation}
\label{2.16.5}
\frac{1}{b-a}\int_a^b f(g(x))dx\geq f\left(\frac{1}{b-a}\int_a^b g(x)dx\right).
\end{equation}

In the absence of this strict convexity, there is no reason to assume homogeneity of the deformation: the deformation could be any function $y(x,t)$ with prescribed values of $y(0,t)=0$ and $y(L,t)$ (given below in (\ref{2.11a})) and with
$\partial y(x,t)/\partial x>1$, where this last condition is required by constancy of the tension $b$ along the rope. The actual deformation which
is selected could depend on higher order gradient terms in the energy function which we have neglected to include, and could be very sensitive to slight inhomogeneities in the rope.

If we assume the selected deformation is homogeneous, then we have
\begin{equation}
\label{2.10}
y(x,t)=\frac{x}{L}y(L,t).
\end{equation}
We can integrate $a(L,t)$ in $t$ to get
\begin{equation}
\label{2.11}
v(L,t)=\int_0^t a(L,t)dt+v_0=a_{0}t+v_0.
\end{equation}
Further integration gives
\begin{equation}
\label{2.11a}
y(L,t)=\int_0^t v(L,t)dt=\frac{a_{0}}{2}t^{2}+v_{0}t+L.
\end{equation}
Thus owing to (\ref{2.10}), we finally arrive at
\begin{equation}
\label{2.12}
y(x,t)=\frac{x}{L}\left(\frac{a_{0}}{2}t^{2}+v_{0}t+L\right).
\end{equation}

\section{The rope with a carabiner}
\label{sec3}

In this section, we consider the climber fall problem in the presence of a carabiner with a friction coefficient $k$.
\begin{figure}[h]
\centering
\includegraphics[width=0.9\textwidth]{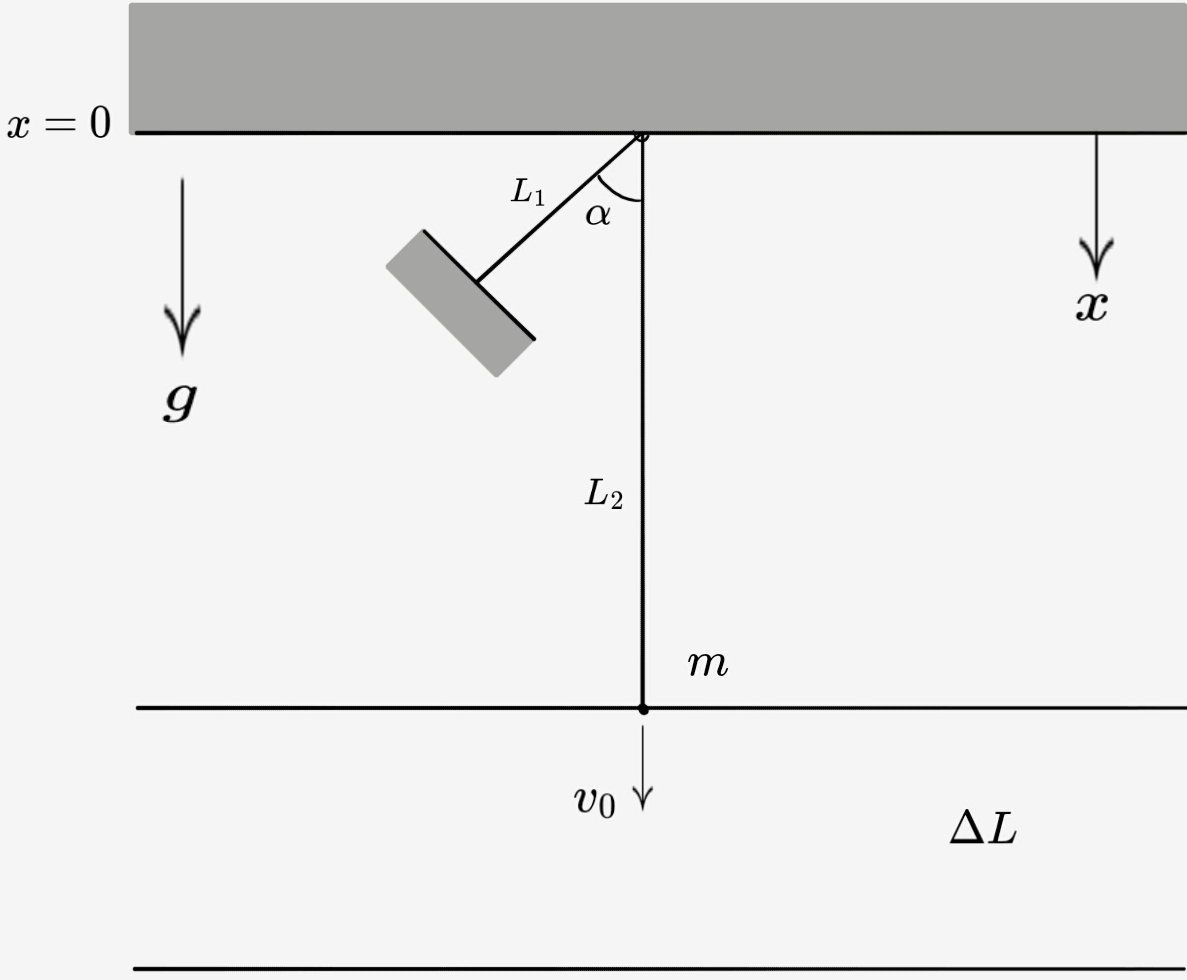}
\caption{\label{Fig3}The configuration after the fall at the instant the rope tightens, in the presence of a carabiner. }
\end{figure}
Namely, we assume that the rope passes through a carabiner attached to the rock wall, with one end of the rope attached to the climber,
while the other end is fastened to the belayer. The rope line from the belayer to the carabiner forms a given angle $\alpha$ with the vertical axis $x,$
see Figure 3. Then we are again seeking a mathematically ideal rope. As was done in Section~\ref{sec2}, we reformulate the problem as follows:\\

\textit{Assume that the rope is vertical and is attached to the rock at the origin $x=0$ and the carabiner is attached to the rock at the point $x=L_1,$ while the rope length between the carabiner and the climber is $L_2$ and the climber has an initial velocity $v_0,$ see Figure 3. In addition, due to the carabiner, the rope tension jumps from $b(t)$ below the carabiner to $\mu b(t)$ above the carabiner, where by the capstan equation (also known as the Euler-Eytelwein formula) $\mu=e^{-(\pi-\alpha)k}$ and $k$ is the coefficient of friction between the metal surface of the carabiner and the rope,}  see Figure 4. \\

\begin{figure}
\centering
\includegraphics[width=0.9\textwidth]{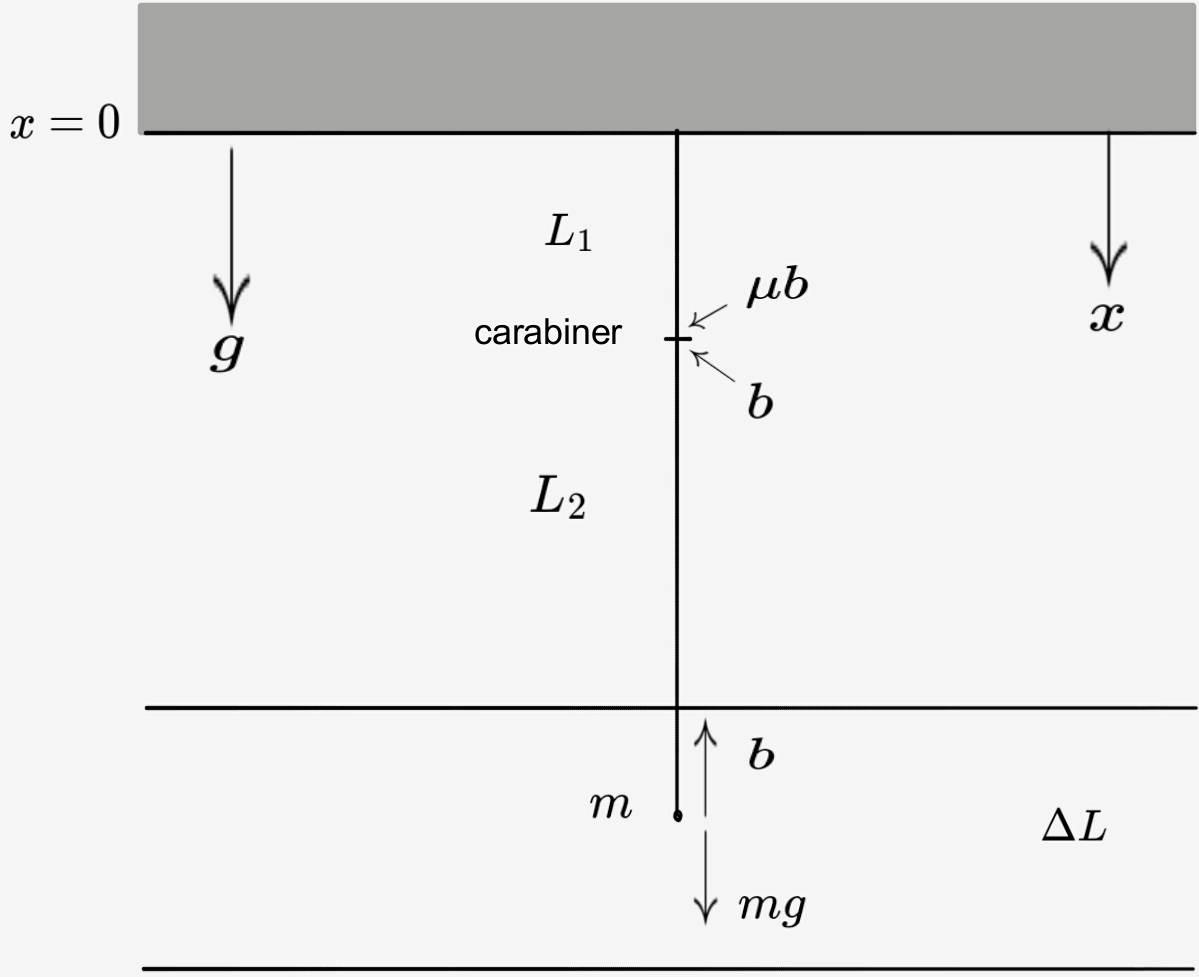}
\caption{\label{Fig4}The reformulated problem for the case with a carabiner.}
\end{figure}
The problem is again finding a rope that solves the problem
\begin{equation}
\label{3.1}
\min(\max_{t\in[0,T]}|b(t)|).
\end{equation}
The setting of the reformulated problem makes it clear that problem~(\ref{3.1}) is solved when $b(t)=b_0$
as well as $a(L,t)=a_0$ for $t\in[0,T]$, where $a_0$ and $b_0$ are given by (\ref{2.8}). We have to now find a rope that develops a constant
resultant force acting on the climber during a fall cycle. It turns out that
the type of rope found in the previous section works also for this case. The idea is to take $L=L_2$ in formula
(\ref{2.16}) and take $y(x,t)=x$ for $x<L_1$ so that the rope develops no deformation between the rock and the carabiner and thus fulfills the desired
conditions. To understand that, let us plot the tension-strain diagram corresponding to the rope given by (\ref{2.16})
with $L$ replaced by $L_2,$ see Figure 5. Because the tension above the carabiner is strictly smaller than
$b_0,$ it is then clear from the diagram that the rope develops no stretching in the segment $[0,L_1]$, as we wanted to establish. Note that that this
conclusion remains valid even if the Euler-Eytelwein formula is questioned, as it may be because the rope diameter is comparable to that of the metal
diameter in the carabiner as observed by Weber and Ehrmann in [\ref{bib:Weber.}]: all we require for the analysis is that the tension in the rope between the carabiner and the belayer be less than that between the carabiner and the climber (i.e. $\mu<1$). The necessity of having the carabiner is to lessen the force on the rope attached to the belayer, to allow the belayer to be positioned in various
places, and to prevent the climber from falling off the wall.
\begin{figure}
\centering
\includegraphics[width=0.9\textwidth]{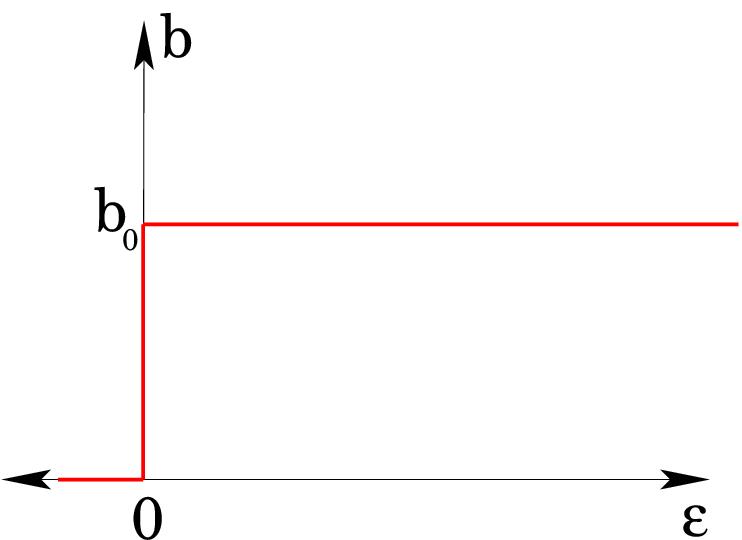}
\caption{\label{Fig5}The tension-strain relation for a mathematically ideal dynamic rope, shown in red. }
\end{figure}

\section{Realizability of the mathematically ideal rope}

\begin{figure}
\centering
\includegraphics[width=0.9\textwidth]{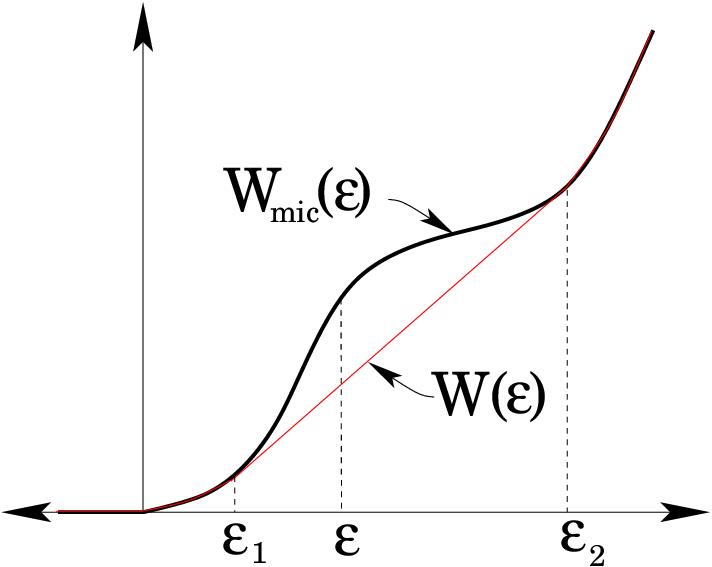}
\caption{\label{Fig6}The black curve is the microscopic energy $W_{mic}$ as a function of the strain $\Gve$, which when convexified gives the macroscopic strain energy, which is the red curve $W(\Gve)$.}
\end{figure}

Here we provide some arguments which suggest that a rope approximately realizing the condition (\ref{2.16})
is not beyond the realm of possibility. The characteristic feature of the tension-strain diagram of Figure 5, is the \textit{plateau
in the tension as the strain is varied.} Wires of shape memory materials such as Nitinol (an alloy of nickel and titanium, see https://en.wikipedia.org/wiki/Nickel\_titanium) have such plateaus, although in currently available shape memory wires the \textit{plateau in the tension} occurs only for strains less than about $8\%$ (see [\ref{bib:Sittner.}] and references therein.) As Fig. 14 in [\ref{bib:Att.}] shows, normal dynamic climbing ropes can have strains of up to $15\%$. The reason for the plateau is illustrated in Figure 6. At a microscopic scale the elastic energy might not be a convex function of $\Gve$ but given by the
function $W_{mic}(\Gve),$ shown in black in Figure 6. A strain having the value $\Gve$ in the region of non-convexity is energetically unstable
and on a macroscopic scale the material \textit{phase separates} into a collection of segments having microscopic strains $\Gve_1$ or $\Gve_2$ in proportions
$\theta$ and $1-\theta$ where $\theta=(\Gve-\Gve_1)/(\Gve_2-\Gve_1)$. The elastic energy density of this mixture is
\begin{eqnarray}
\label{6.0a}
W(\Gve) & = & \theta W_{mic}(\Gve_1)+(1-\theta)W_{mic}(\Gve_2) \\ \nonumber
        & = & W_{mic}(\Gve_2)+\frac{\Gve-\Gve_1}{\Gve_2-\Gve_1}[W_{mic}(\Gve_1)-W_{mic}(\Gve_2)],
\end{eqnarray}
which depends linearly on $\Gve$ for strains between $\Gve_1$ and $\Gve_2$: thus the wire tension is independent of the
macroscopic strain in this interval. (We remark in passing that the treatment for the deformation of two or three dimensional shape memory materials,
rather than wires, has also been developed but is more complicated and involves quasiconvexification rather than convexification: see, for example, [\ref{bib:Ball.}].)

\begin{figure}
\centering
\includegraphics[width=0.9\textwidth]{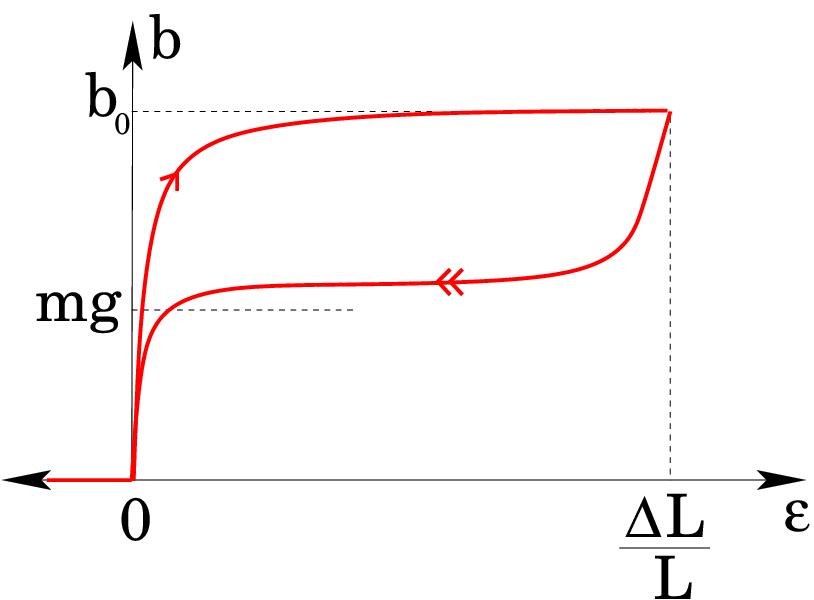}
\caption{\label{Fig7}The tension-strain hysteresis loop for a close to ideal dynamic rope, shown in red. The path with the single red arrow denotes the trajectory when the strain is increased, the double red arrow denotes the trajectory when the strain is decreased.}
\end{figure}

In reality, the tension-strain diagram of shape memory wires have a hysteresis loop, which is not
encompassed by our purely elastic formulation. Although hysteresis can be minimized in shape memory materials [\ref{bib:Song.}], a
hysteresis loop in the response of a rope would actually be a necessity: it would be good if the tension in the rope was
just slightly more than $mg$ when the rope retracts after it reaches its maximum extension, as then the climber would slowly
rebound from the fall, returning close to $x=L$. This close-to-ideal tension-strain relation with a hysteresis loop is sketched in Figure 7.
In some circumstances it might be best if the rope remained close to the maximum extension $x=L+\Delta L$ after the fall: this will be the
case if the return path in the hysteresis loop was below the line $b=mg$.

\section{Conclusions}
\label{sec3}

We do not expect this paper to have an immediate effect on the climbing community, but by providing a prescription for a \textit{mathematically ideal rope,}
the work may help guide the development of new ropes. Also, it may motivate research into alternate shape memory materials which may
be more suitable for use in ropes than nitinol. It is worth
mentioning, that there are recently developed shape memory materials that
are not alloys, but polymers, e.g., [\ref{bib:Eis.Rhe.},\ref{bib:Xie.Rou.}], although these do not yet have the required strength to be
suitable in ropes.
The suggestion of using materials with shape memory wire characteristics appears to be novel.
However, in reality the tension-strain hysteresis loop trajectory in shape memory materials depends on the rate of
deformation, and the plateaus disappear if the rate of deformation is too fast (see, for example, figure 14 in [\ref{bib:Heller.}]).
In order for the plateaus to be retained at a desired rate of deformation, the \textit{phase separation} must occur on an appropriately fast time scale. Another unwanted characteristic is that the shape
and position of the hysteresis loop in general varies according to the number of deformation cycles the wire has undergone.

While shape memory wires only have stress plateau extending up to 8\% strain, while dynamic ropes can have strains of 15\%, this may actually be an
advantage as pointed out to us by Wendy Crone at the University of Wisconsin. The reason is that with a \textit{mathematically ideal rope,} one
can have the same peak force with less overall elongation than in a standard dynamic rope, and less elongation minimizes the chance of collision
with a rock outcrop or another climber.

Of course a ``rope'' built from
shape memory material wires may have disadvantages not considered here, such as being too heavy (although titanium is comparatively light among metals),
too expensive, not easily coiled, not easily knotted, or having properties which are temperature-dependent. Cables have been built from shape memory
wires and their properties have been studied [\ref{bib:Reedlunn.}, \ref{bib:Ozbulut.}]. Significantly, Reedlunn, Daly, and Shaw [\ref{bib:Reedlunn.}]
find that the hysteresis loop can extend
up to a strain of about 12\% (see their Figure 6a), depending on the placement of the wires in the cable, but in this case the stress plateau is lost:
an alternate cable design has (as shown in the same figure) a hysteresis loop up to a strain of about 8\% while retaining the stress plateau. Cables, like
wires, have a hysteresis loop that varies according to the number of deformation cycles: see Figure 9 of Ozbulut, Daghash,
and Sherif [\ref{bib:Ozbulut.}] where the hysteresis loop after
100 cycles is much narrower than after one cycle. It may be an advantage to combine
fibers or wires with shape memory characteristics with other rope elements, but we have not studied this possibility.

We remark that recently developed shape memory alloys can reliably
go through 10 million transformation cycles without fatigue [\ref{bib:SMfat}]; although in this case the hysteresis loop is quite small.
There may be other applications where this analysis is useful. For example, suppose one wants to drop cargo from a plane or a helicopter,
without deploying parachutes (which may easily be seen by an enemy).
Then, to lessen the impact of the fall on the cargo, it may be useful to attach a tether to the cargo which becomes taut
at a certain distance above the ground: the distance to the ground sets the maximal elongation.
One would again want to minimize the peak tension in the tether, thus minimizing the forces on the plane or the helicopter
and on the cargo. In such applications to achieve the desired \textit{mathematically ideal tether} which produces a constant braking force,
one could use a tether with sacrifical elements which tear and dissipate energy as it is stretched. A one dimensional model with these
sacrificial elements was studied by Cherkaev, Cherkaev and Slepyan [\ref{bib:Cherkaev}] and shows an approximate plateau in the tension verses
elongation (with oscillations): see their figure 3. Such sacrificial elements are believed to be a mechanism giving bone its strength
[\ref{bib:Thompson.}] and account for resilance of double network hydrogels [\ref{bib:Gong}: see also http://imechanica.org/node/13088].

\section*{Acknowledgements}

The authors are grateful to Richard James for drawing their attention to fatigue resistant shape memory material alloys, to Tom Shield for drawing
their attention to work on shape memory material cables, and to Hong Wei for drawing their attention to double network hydrogels.
Additionally, Wendy Crone is thanked for pointing out that the shorter elongation of 8\% (within the stress plateau) of shape memory
material wires can be an advantage if the peak force is less than dynamic ropes that have a 15\% elongation. The referees are thanked for their
many useful and insightful comments on the manuscript which led to substantial improvements.

\appendix{Appendix I}\\

\textit{\textbf{Notation}}

\hspace*{-0.3cm}$a_0$\hspace{1.5cm} constant acceleration/deceleration\\
\hspace*{2.0cm}of the climber\\
$a(x,t)$\hspace{0.9cm} acceleration of the rope point at time $t$\\
\hspace*{1.9cm} with initial (at $t=0$) coordinate $x$\\
$b(t)$\hspace{1.3cm} rope tension at the point $x=L$ at\\
\hspace*{1.8cm} time $t$\\
$b_0$\hspace{1.5cm} constant rope tension at the point $x=L$\\
\hspace*{1.8cm} at constant acceleration\\
$E_{el}$ \hspace{1.2cm} stored elastic energy of the rope\\
$E_{total}$ \hspace{1.0cm} total energy of the system\\
$\epsilon(x), \epsilon_i$ \hspace{0.8cm} elastic strain\\
$g$ \hspace{1.7cm} gravitational acceleration\\
$h_0$ \hspace{1.5cm} initial coordinate of the climber\\
$k$ \hspace{1.6cm} friction coefficient between the rope\\
\hspace*{2.0cm} and the carabineer\\
$L$ \hspace{1.6cm} length of the undeformed rope\\
$\triangle L$ \hspace{1.3cm} maximal stretch of the rope in a\\
\hspace*{2.0cm} fall cycle\\
$m$ \hspace{1.6cm} mass of the climber\\
$t$ \hspace{1.7cm} time\\
$T$ \hspace{1.6cm} time moment when the rope reaches\\
\hspace*{2.0cm}  maximal stretch for the first time\\
$u(x,t)$ \hspace{0.9cm} the rope displacement function\\
$v(x,t)$ \hspace{1.0cm}velocity of the rope point $x$ ate time $t$\\
$v_0$ \hspace{1.5cm} initial velocity of the climber in the\\
\hspace*{2.0cm} reformulated problem\\
$W$ \hspace{1.6cm} rope elastic energy density\\
$x$ \hspace{1.6cm} coordinate variable of the system\\
$y(x,t)$ \hspace{0.9cm} the rope deformation function\\

\end{document}